\documentclass[11pt]{article}

\usepackage[T1]{fontenc}
\usepackage{lmodern}
\usepackage[margin=1.05in]{geometry}
\usepackage{needspace}
\usepackage{amsmath,amssymb,amsthm,mathtools}
\usepackage{mathrsfs}
\usepackage[colorlinks=true,linkcolor=blue,citecolor=blue,urlcolor=blue]{hyperref}

\numberwithin{equation}{section}

\newcommand{\R}{\mathbb{R}}
\newcommand{\eps}{\varepsilon}
\renewcommand{\Cap}{\operatorname{Cap}}
\newcommand{\dist}{\operatorname{dist}}
\newcommand{\supp}{\operatorname{supp}}
\newcommand{\loc}{\mathrm{loc}}
\newcommand{\dd}{\,d}
\newcommand{\one}{\mathbf{1}}

\newcommand{\cD}{\mathcal D}
\newcommand{\cL}{\mathcal L}
\newcommand{\cA}{\mathcal A}
\newcommand{\cE}{\mathcal E}
\newcommand{\cH}{\mathcal H}
\newcommand{\sstar}{2_s^*}
\newcommand{\hstar}{2_{\mu,s}^*}

\newtheorem{theorem}{Theorem}[section]
\newtheorem{lemma}[theorem]{Lemma}
\newtheorem{proposition}[theorem]{Proposition}
\newtheorem{corollary}[theorem]{Corollary}

\theoremstyle{definition}
\newtheorem{definition}[theorem]{Definition}

\theoremstyle{remark}
\newtheorem{remark}[theorem]{Remark}

\title{On symmetry and monotonicity of singular solutions to \\critical fractional Hartree equations}

\author{Guangze Gu, Aleks Jevnikar, Zhipeng Yang\thanks{Corresponding author: yangzhipeng326@163.com}}

\date{}

\AtEndDocument{%
\par
\bigskip
\bigskip
\noindent\begin{minipage}{\textwidth}

\noindent
\textbf{Guangze Gu}\\[0.2em]
\textsc{Department of Mathematics, Yunnan Normal University, Kunming 650500, China}\\
\textsc{Yunnan Key Laboratory of Modern Analytical Mathematics and Applications, Kunming 650500, China}\\[0.3em]
\textit{E-mail address}: \texttt{guangzegu@163.com}

\bigskip

\noindent
\textbf{Aleks Jevnikar}\\[0.2em]
\textsc{Department of Mathematics, Computer Science and Physics, University of Udine, Udine 33100, Italy}\\[0.3em]
\textit{E-mail address}: \texttt{aleks.jevnikar@uniud.it}

\bigskip

\noindent
\textbf{Zhipeng Yang}\\[0.2em]
\textsc{Department of Mathematics, Yunnan Normal University, Kunming 650500, China}\\
\textsc{Yunnan Key Laboratory of Modern Analytical Mathematics and Applications, Kunming 650500, China}\\[0.3em]
\textit{E-mail address}: \texttt{yangzhipeng326@163.com}
\end{minipage}
}

\begin{document}

\maketitle

\begin{abstract}
This paper continues the study initiated in \cite{GuJevnikarYang2026b} of fractional Hartree equations with zero capacity singular sets. We are concerned here with the critical case, i.e.
\[
(-\Delta)^s u=\left(\int_{\R^N\setminus\Gamma}
\frac{u^{\hstar}(y)}{|x-y|^\mu}\dd y\right)u^{\hstar-1}
\quad\text{in }\R^N\setminus\Gamma,
\]
where $N\ge2$, $s\in(0,1)$, $\mu\in(0,N)$ and $\hstar$ is the critical exponent associated to this class of problems. We prove symmetry and monotonicity with respect to the singular set $\Gamma$.  The Kelvin transforms centered at regular points provide sufficient decay at infinity for both the transformed solutions and their associated Hartree coefficients, without imposing global integrability assumptions. Combined with a moving plane argument and a comparison estimate in critical Lebesgue spaces, this allows the moving plane method to cover the full interaction range.
\end{abstract}



\section{Introduction and Main Results}\label{Sec1}

Let $N\ge2$, $s\in(0,1)$, and $\mu\in(0,N)$, and set
\[
\sstar=\frac{2N}{N-2s}, \qquad \hstar=\frac{2N-\mu}{N-2s}.
\]
We consider positive weak singular solutions of
\begin{equation}\label{eq1.1}
(-\Delta)^s u=\cH[u]u^{\hstar-1} \quad\text{in }\R^N\setminus\Gamma, \qquad \cH[u](x)=\int_{\R^N\setminus\Gamma}\frac{u^{\hstar}(y)}{|x-y|^\mu}\dd y,
\end{equation}
where $\Gamma$ is compact, contained in the hyperplane $\{x_1=0\}$, and has zero $(s,2)$-capacity. We prove hyperplane symmetry and strict transverse monotonicity throughout the interaction range $0<\mu<N$, including the case $\mu>4s$ in which the outer exponent $\hstar-1$ lies below one.

The exponent $\hstar$ is determined by the Hardy--Littlewood--Sobolev inequality and the scaling
\[
u_\rho(x)=\rho^{\frac{N-2s}{2}}u(\rho x).
\]
This scaling leaves invariant both the homogeneous fractional energy and the Hartree integral
\[
\iint_{\R^N\times\R^N} \frac{|u(x)|^{\hstar}|u(y)|^{\hstar}}{|x-y|^\mu}\dd x\dd y.
\]

The critical scaling also identifies the singular profile $|x|^{-(N-2s)/2}$, which solves the equation after multiplication by a suitable positive constant; see Remark~\ref{Rem2.4}. Its $\sstar$-power is proportional to $|x|^{-N}$ and is therefore nonintegrable both at the origin and at infinity. A symmetry theorem for singular critical solutions must consequently accommodate functions outside the global critical Lebesgue space. For a general singular set, a further difficulty is that reflection moves a copy of that set into the comparison region. Zero fractional capacity allows this reflected set to be removed by energy cutoffs, but it does not by itself control the Hartree interaction or the behavior at infinity. The result below addresses these three issues within the weak formulation.

The classical Choquard problem studied by Lieb \cite{Lieb1977} has $N=3$, the Newton kernel, and $p=2$; Lions \cite{Lions1980} developed its variational theory further. For the general family
\[
-\Delta u+u=(I_\alpha*|u|^p)|u|^{p-2}u \quad\text{in }\R^N,
\]
where $I_\alpha$ is the Riesz potential, Moroz and Van Schaftingen \cite{MorozVanSchaftingen2013} established existence, qualitative properties, and decay of ground states. They also studied the lower HLS critical exponent with an external potential \cite{MorozVanSchaftingen2015} and surveyed the theory in \cite{MorozVanSchaftingen2017}. For fractional Choquard equations, d'Avenia, Siciliano, and Squassina \cite{dAveniaSicilianoSquassina2015} developed the ground-state theory. Gao and Yang \cite{GaoYang2017} and Mukherjee and Sreenadh \cite{MukherjeeSreenadh2017} treated critical problems for the local and fractional operators, respectively. These variational problems rely on the sharp Hardy--Littlewood--Sobolev inequality \cite{Lieb1983,LiebLoss2001}.

Gidas, Ni, and Nirenberg \cite{GidasNiNirenberg1979} applied moving planes to semilinear elliptic equations, and Berestycki and Nirenberg \cite{BerestyckiNirenberg1991} developed moving plane and sliding arguments for further classes of equations. Li \cite{Li2004} and Chen, Li, and Ou \cite{ChenLiOu2006} used moving spheres and moving planes in integral form for conformally invariant equations. Chen, Li, and Li \cite{ChenLiLi2017} developed a direct method of moving planes for the fractional Laplacian, and Chen, Li, and Zhang \cite{ChenLiZhang2017} developed a direct method of moving spheres. For fractional Hartree equations, Dai, Fang, and Qin \cite{DaiFangQin2018} used direct moving planes, Dai, Huang, Qin, Wang, and Fang \cite{DaiHuangEtAl2019} used integral moving planes, and Dai, Liu, and Qin \cite{DaiLiuQin2021} used integral moving spheres. These whole-space classification arguments must be modified when the equation holds only outside $\Gamma$, since reflection places a copy of the singular set inside the moving half-space.

Brezis and V\'eron \cite{BrezisVeron1980} studied removable singularities for local semilinear equations, while Caffarelli, Gidas, and Spruck \cite{CaffarelliGidasSpruck1989} analyzed the critical isolated-singularity problem. Caffarelli, Jin, Sire, and Xiong \cite{CaffarelliJinSireXiong2014} developed the corresponding local analysis for the fractional critical equation. For compact singular sets of zero capacity, Montoro, Punzo, and Sciunzi \cite{MontoroPunzoSciunzi2018} proved symmetry and monotonicity for fractional equations with a local nonlinearity $f(u)$. Their whole-space theorem assumes global $L^1$ integrability and exterior $L^{\sstar}$ integrability. More recently, Du and Yang \cite{DuYang2026} proved hyperplane symmetry and transverse monotonicity for positive classical solutions of $-\Delta u=f(u)$ in dimension $N\ge3$, with a nonremovable singularity on a closed set of zero Newtonian capacity contained in the hyperplane. They assume that $f$ is locally bounded on $(0,\infty)$ and $t^{-(N+2)/(N-2)}f(t)$ is nonincreasing, and use moving spheres without an assumption at infinity. The Hartree term in \eqref{eq1.1} requires an additional comparison for the reflected convolution potential.

For local Choquard equations, Chen and Zhou \cite{ChenZhou2016} classified isolated singularities; in \cite{ChenZhou2023} they also studied qualitative properties of singular solutions. Andrade, Feng, Piccione, and Yang \cite{AndradeFengPiccioneYang2026a} proved that positive classical solutions of the critical local Hartree equation in the punctured space that tend to infinity at the origin are radially symmetric and radially nonincreasing. For the punctured-ball problem, they obtained a pointwise upper bound and asymptotic radial symmetry. Cai, Gu, and Jevnikar \cite{CaiGuJevnikar2025} treated a class of local Hartree equations with convex nonlinearities, for $N\ge3$ and $N-2<\mu<N$, and proved hyperplane symmetry and monotonicity when the singular set is a compact subset of a hyperplane with zero $2$-capacity. 

For the fractional equation with independent powers, Gu, Jevnikar, and Yang \cite{GuJevnikarYang2026a} derived a Riesz decomposition that retains a nonnegative fundamental-solution term. Their symmetry theorem assumes $p,q\ge1$, $p\le\hstar$, $q\le\hstar-1$, an upper bound of fundamental-solution order at the origin, weighted integrability of the source, and $u(x)\to+\infty$ as $x\to0$. An off-center moving-spheres argument then yields radial symmetry and strict radial decrease. For $p=\hstar$ and $q=\hstar-1$, Andrade, Feng, Piccione, and Yang \cite{AndradeFengPiccioneYang2026b} proved radial symmetry for positive $C^{1,1}_{\loc}(\R^N\setminus\{0\})\cap\cL_s(\R^N)$ solutions with a nonremovable singularity, for every $0<\mu<N$. Their proof uses the Caffarelli--Silvestre extension and moving spheres and also yields radial nonincrease.

Gu, Jevnikar, and Yang \cite{GuJevnikarYang2026b} treated compact zero-capacity singular sets for the subcritical equation
\[
(-\Delta)^s u=(|\cdot|^{-\mu}*u^p)u^q.
\]
Their Theorems~1.1--1.2 assume $p,q>1$, $p+q+1<2\hstar$, and suitable global Lebesgue integrability. Their Theorem~1.3 instead assumes
\[
1<p\le\hstar,\qquad 1<q\le\hstar-1,
\qquad (p,q)\ne(\hstar,\hstar-1),
\]
and obtains symmetry from the weak solution class without an additional global Lebesgue condition. Local H\"older continuity is derived in both settings. The simultaneous critical equality is excluded from these results; moreover, $q>1$ in the second set of inequalities forces $\mu<4s$.

The present paper treats precisely the conformally invariant pair $p=\hstar$, $q=\hstar-1$ for continuous weak solutions and for every $0<\mu<N$. Theorem~\ref{Thm1.2} allows any compact set of zero $(s,2)$-capacity contained in the reflecting hyperplane, with a single nonremovable point sufficient to fix the symmetry plane. Relative to the isolated-singularity result of \cite{AndradeFengPiccioneYang2026b}, the principal extension is the singular-set geometry and the weak formulation. Relative to \cite{GuJevnikarYang2026b}, the theorem includes the critical equality pair and the range $\mu\ge4s$ when that range is nonempty. 

We now state the weak formulation and the main theorem. For $v:\R^N\to\R$, write
\[
[v]_{H^s(\R^N)}^2=\iint_{\R^N\times\R^N}\frac{|v(x)-v(y)|^2}{|x-y|^{N+2s}}\dd x\dd y,
\qquad \|v\|_{H^s(\R^N)}^2=\|v\|_{L^2(\R^N)}^2+[v]_{H^s(\R^N)}^2,
\]
and let $H^s(\R^N)$ be the space of functions with finite norm. For an open set $\Omega$, the notation $v\in H^s_{\loc}(\Omega)$ means that $\chi v\in H^s(\R^N)$ for every $\chi\in C_c^\infty(\Omega)$, with $\chi v$ extended by zero outside $\Omega$. For a compact set $K\subset\R^N$, let
\[
\Cap_{s,2}(K)=\inf\left\{ \|\varphi\|_{H^s(\R^N)}^2: \varphi\in C_c^\infty(\R^N),\quad \varphi\ge1\text{ in a neighborhood of }K \right\},
\]
and set
\[
\cL_s(\R^N)=\left\{ v\in L^1_{\loc}(\R^N): \int_{\R^N}\frac{|v(x)|}{1+|x|^{N+2s}}\dd x<\infty \right\}.
\]
Let $C_{N,s}>0$ denote the standard normalization constant for $(-\Delta)^s$. For $\phi\in C_c^\infty(\R^N\setminus\Gamma)$, write
\[
\cE(v,\phi)=\frac{C_{N,s}}2 \iint_{\R^N\times\R^N} \frac{(v(x)-v(y))(\phi(x)-\phi(y))}{|x-y|^{N+2s}}\dd x\dd y.
\]

\begin{definition}\label{Def1.1}
A measurable function $u$ is a positive weak singular solution of \eqref{eq1.1} if
\[
u>0\quad\text{almost everywhere in }\R^N\setminus\Gamma, \qquad u\in\cL_s(\R^N)\cap H^s_{\loc}(\R^N\setminus\Gamma)\cap C(\R^N\setminus\Gamma),
\]
$\cH[u](x)<\infty$ for almost every $x\in\R^N\setminus\Gamma$, $\cH[u]u^{\hstar-1}\in L^1_{\loc}(\R^N\setminus\Gamma)$, and
\[
\cE(u,\phi)=\int_{\R^N\setminus\Gamma}\cH[u]u^{\hstar-1}\phi\dd x
\]
for every $\phi\in C_c^\infty(\R^N\setminus\Gamma)$. A point $x_0\in\Gamma$ is nonremovable if
\[
\limsup_{\R^N\setminus\Gamma\ni x\to x_0}u(x)=+\infty.
\]
Points of $\R^N\setminus\Gamma$ are called regular points of $u$.
\end{definition}

\begin{theorem}\label{Thm1.2}
Let $N\ge2$, $s\in(0,1)$, and $\mu\in(0,N)$. Let $\Gamma\subset\{x_1=0\}$ be compact with $\Cap_{s,2}(\Gamma)=0$. Suppose that $u$ is a positive weak singular solution of \eqref{eq1.1} in the sense of Definition~\ref{Def1.1} and has a nonremovable singular point on $\Gamma$. Then
\[
u(x_1,x')=u(-x_1,x') \quad\text{for every }(x_1,x')\in\R^N\setminus\Gamma.
\]
Moreover, for every $x'\in\R^{N-1}$,
\[
u(x_1,x')<u(y_1,x') \quad\text{whenever }x_1<y_1<0,
\]
and
\[
u(y_1,x')<u(x_1,x') \quad\text{whenever }0<x_1<y_1.
\]
\end{theorem}

\begin{corollary}\label{Cor1.3}
Under the assumptions of Theorem~\ref{Thm1.2}, if $\Gamma=\{0\}$, then $u$ is radially symmetric about the origin. More precisely, there exists a strictly decreasing function $U:(0,\infty)\to(0,\infty)$ such that
\[
u(x)=U(|x|) \quad\text{for every }x\in\R^N\setminus\{0\}.
\]
\end{corollary}

Since no decay condition is imposed on $u$, we begin with a Kelvin transform centered at a regular point $a\in\{x_1=0\}\setminus\Gamma$:
\[
\mathcal I_a(x)=a+\frac{x-a}{|x-a|^2}, \qquad u_a(x)=|x-a|^{-(N-2s)}u(\mathcal I_a(x)).
\]
The critical equation is invariant under this transformation, and the transformed singular set
\[
\Gamma_a=\mathcal I_a(\Gamma)\cup\{a\}
\]
remains in the same hyperplane and has zero $(s,2)$-capacity. Continuity at $a$ gives $u_a(x)=O(|x|^{-(N-2s)})$ at infinity. To cover the full interaction range, the comparison argument also needs a tail estimate for
\[
\cA_u=\cH[u]u^{\hstar-2},
\]
initially defined almost everywhere, with $\cA_u=\cH[u]$ when $\hstar=2$. We first derive strict positivity at every regular point from the weak equation. Thus $\cA_u$ is locally bounded even when $\hstar-2<0$. The critical exponent then gives the exact covariance
\[
\cA_{u_a}(x)=|x-a|^{-4s}\cA_u(\mathcal I_a(x)),
\qquad \mu+(N-2s)(\hstar-2)=4s.
\]
Consequently $\cA_{u_a}(x)=O(|x|^{-4s})$ at infinity. Both $u_a^{\sstar}$ and $\cA_{u_a}^{N/(2s)}$ are integrable there, with vanishing exterior integrals. These two tail properties are proved for each Kelvin transform; neither is an assumption on the original solution.

The second ingredient is the scalar estimate
\[
t^q-r^q\le\max\{1,q\}\,t^{q-1}(t-r),
\qquad q>0,\quad 0\le r<t,
\]
proved in Lemma~\ref{Lem2.8}. When $0<q\le1$, it follows by dividing by $t^q$ and using $(r/t)^q\ge r/t$. Applied with $q=\hstar-1$, it keeps the coefficient $\cA_u$ at the larger comparison value for both the convex and concave powers. This permits one comparison argument throughout $0<\mu<N$.

More precisely, Proposition~\ref{Prop3.5} requires smallness, on the set where the solution exceeds its reflection, of the $L^{\sstar}$-norm of the solution and the $L^{N/(2s)}$-norm of its coefficient. We subtract the reflected equation and test with the positive part of the difference, multiplied by cutoffs near the reflected singular set and at infinity. An antisymmetric energy inequality controls the fractional term, while the Hardy--Littlewood--Sobolev and H\"older inequalities control the Hartree term. Fatou's lemma then permits removal of the cutoffs. Applied to $u_a$, this estimate starts the planes far from the singular set and, together with strict comparison and the nonremovable singularity, allows them to reach the hyperplane.

Symmetry passes from $u_a$ to $u$ because inversion centered on the hyperplane commutes with reflection. Transverse monotonicity requires a separate argument: inversion maps a transverse line not passing through $a$ onto a circle through $a$, so it does not preserve the order along such lines. We let the center $a$ tend to infinity within the hyperplane and recover weak monotonicity of $u$; the strict comparison principle, which requires neither critical tail, then gives strict monotonicity. Since $N\ge2$ and $\Gamma$ is compact, the hyperplane contains arbitrarily distant regular centers. The same dimensional assumption implies $N>2s$, which ensures that the added point $\{a\}$ has zero capacity.

Section~\ref{Sec2} contains the capacity, Kelvin transform, and integrability estimates. Section~\ref{Sec3} proves the reflected comparison estimate. Section~\ref{Sec4} carries out the moving plane argument under the two critical tail conditions, and Section~\ref{Sec5} transfers symmetry and monotonicity from the Kelvin transforms to the original solution.

\section{Kelvin Transforms and Critical Estimates}\label{Sec2}

For an open set $\Omega\subset\R^N$, we write
\[
H_0^s(\Omega)= \{v\in H^s(\R^N):v=0\text{ almost everywhere in }\R^N\setminus\Omega\}.
\]
We will use the following localization fact. If $v\in\cL_s(\R^N)\cap H^s_{\loc}(U)$, $K\Subset U$, and $\phi\in H^s(\R^N)$ is supported in $K$, then $\cE(v,\phi)$ is absolutely convergent and depends continuously on $\phi$ in $H^s$ with its support in a fixed regular compact set. Indeed, choose $\chi\in C_c^\infty(U)$ equal to one near $K$. The pairing with $\chi v\in H^s$ is controlled by Cauchy--Schwarz. The remaining contribution equals
\[
-C_{N,s}\int_K\phi(x)\int_{\R^N}
\frac{(1-\chi(y))v(y)}{|x-y|^{N+2s}}\dd y\dd x.
\]
The inner integral is uniformly bounded in absolute value on $K$: separation controls the bounded region using $v\in L^1_{\loc}$, and the $\cL_s$ condition controls infinity. This term is therefore bounded by $C_K\|\phi\|_{L^1(K)}$. If the weak source is locally bounded, the weak equation consequently extends to these tests by mollification. Nonnegative tests admit nonnegative mollifications supported in a slightly larger fixed compact subset of $U$. We will use this fact only with supports strictly inside the regular set.

\begin{lemma}\label{Lem2.1}
Let $K\subset\R^N$ be compact and satisfy $\Cap_{s,2}(K)=0$. For every open neighborhood ${O}$ of $K$ and every $\eps>0$, there exists $\eta_\eps\in C^\infty(\R^N)$ such that
\[
0\le\eta_\eps\le1, \qquad \eta_\eps=0\text{ in a neighborhood of }K, \qquad \eta_\eps=1\text{ on }\R^N\setminus{O},
\]
and
\[
\|1-\eta_\eps\|_{H^s(\R^N)}<\eps.
\]
Moreover, $|K|=0$, every Euclidean reflection of $K$ has zero $(s,2)$-capacity, and a finite union of compact sets of zero capacity has zero capacity.
\end{lemma}

\begin{proof}
Choose open sets $K\Subset{O}_1\Subset{O}$ and $\chi\in C_c^\infty({O})$ with $0\le\chi\le1$ and $\chi=1$ on ${O}_1$. By the definition of capacity, there are $\theta_j\in C_c^\infty(\R^N)$ such that $\theta_j\ge1$ near $K$ and $\|\theta_j\|_{H^s}\to0$. Since multiplication by $\chi$ is bounded on $H^s(\R^N)$, the functions $\psi_j=\chi\theta_j$ satisfy $\|\psi_j\|_{H^s}\to0$, are supported in ${O}$, and are at least one near $K$. Let $\Theta\in C^\infty(\R)$ be Lipschitz, with
\[
0\le\Theta\le1, \qquad \Theta(t)=0\text{ for }t\le\frac12, \qquad \Theta(t)=1\text{ for }t\ge1.
\]
The Lipschitz bound and $\Theta(0)=0$ give $\|\Theta(\psi_j)\|_{H^s}\le C\|\psi_j\|_{H^s}$. Thus $\eta_j=1-\Theta(\psi_j)$ has the required properties for large $j$.

For every admissible function $\varphi$ in the definition of capacity,
\[
|K|\le\|\varphi\|_{L^2(\R^N)}^2\le\|\varphi\|_{H^s(\R^N)}^2.
\]
The infimum of the right-hand side is zero, so $|K|=0$. Reflections preserve Lebesgue measure, distances, and the $H^s$-norm, hence preserve capacity. If $K_1$ and $K_2$ have zero capacity, use the preceding construction to choose nonnegative functions $\varphi_i\in C_c^\infty(\R^N)$ that equal one near $K_i$ and have arbitrarily small $H^s$-norms. The sum $\varphi_1+\varphi_2$ is admissible for $K_1\cup K_2$, and
\[
\|\varphi_1+\varphi_2\|_{H^s}\le\|\varphi_1\|_{H^s}+\|\varphi_2\|_{H^s}.
\]
Thus $K_1\cup K_2$ has zero capacity. The assertion for finite unions follows by induction. For background on capacity, see \cite{AdamsHedberg1996}.
\end{proof}

For $a\in\R^N$, let
\[
\mathcal I_a(x)=a+\frac{x-a}{|x-a|^2}, \qquad x\ne a.
\]
The map $\mathcal I_a$ is an involution and satisfies
\begin{equation}\label{eq2.1}
|\mathcal I_a(x)-\mathcal I_a(y)|=\frac{|x-y|}{|x-a||y-a|}, \qquad |\det D\mathcal I_a(x)|=|x-a|^{-2N}.
\end{equation}

\begin{lemma}\label{Lem2.2}
Let $K\subset\R^N$ be compact, $a\notin K$, and $\Cap_{s,2}(K)=0$. Then
\[
\Cap_{s,2}\bigl(\mathcal I_a(K)\cup\{a\}\bigr)=0.
\]
If $K\subset\{x_1=0\}$ and $a\in\{x_1=0\}\setminus K$, then
\[
\mathcal I_a(K)\cup\{a\}\subset\{x_1=0\}.
\]
\end{lemma}

\begin{proof}
Choose bounded open sets $K\Subset{O}_0\Subset{O}\Subset\R^N\setminus\{a\}$. Lemma~\ref{Lem2.1} provides admissible functions $\psi_j\in C_c^\infty({O}_0)$ with $\|\psi_j\|_{H^s}\to0$. Extend $v_j=\psi_j\circ\mathcal I_a$ by zero at $a$ and put $D=\mathcal I_a({O})$. Then $v_j\in C_c^\infty(D)$, is admissible for $\mathcal I_a(K)$, and has support in the fixed compact set $\mathcal I_a(\overline{{O}_0})\Subset D$. The Jacobians and distance ratios in \eqref{eq2.1} are bounded above and below on these neighborhoods, so
\[
\|v_j\|_{L^2}^2+\iint_{D\times D}\frac{|v_j(x)-v_j(y)|^2}{|x-y|^{N+2s}}\dd x\dd y\le C\|\psi_j\|_{H^s}^2.
\]
The remaining pairs satisfy
\[
\iint_{D\times D^c}\frac{|v_j(x)|^2}{|x-y|^{N+2s}}\dd x\dd y\le C\|v_j\|_{L^2}^2,
\]
because the supports have a fixed positive distance from $D^c$. Thus $\|v_j\|_{H^s}\to0$, proving that $\mathcal I_a(K)$ has zero capacity.

The standing assumptions $N\ge2$ and $s\in(0,1)$ imply $N>2s$, so a singleton has zero $(s,2)$-capacity. Indeed, if $\zeta\in C_c^\infty(B_2)$ equals one on $B_1$, then $\zeta_\rho(x)=\zeta((x-a)/\rho)$ satisfies
\[
\|\zeta_\rho\|_{L^2}^2=O(\rho^N), \qquad [\zeta_\rho]_{H^s}^2=O(\rho^{N-2s}).
\]
The union therefore has zero capacity by Lemma~\ref{Lem2.1}. The last assertion follows directly from the formula for $\mathcal I_a$, since the first coordinate remains zero on the hyperplane.
\end{proof}

\begin{proposition}\label{Prop2.3}
Let $u$ be a positive weak singular solution of \eqref{eq1.1}, let $a\in\{x_1=0\}\setminus\Gamma$, and set
\[
\Gamma_a=\mathcal I_a(\Gamma)\cup\{a\}, \qquad u_a(x)=|x-a|^{-(N-2s)}u(\mathcal I_a(x)), \quad x\in\R^N\setminus\Gamma_a.
\]
Then $\Gamma_a$ is a compact subset of $\{x_1=0\}$ with zero $(s,2)$-capacity, and $u_a$ is a positive weak singular solution of the same critical equation,
\begin{equation}\label{eq2.2}
(-\Delta)^s u_a=\cH[u_a]u_a^{\hstar-1} \quad\text{in }\R^N\setminus\Gamma_a.
\end{equation}
Moreover,
\begin{align}
\cH[u_a](x)&=|x-a|^{-\mu}\cH[u](\mathcal I_a(x)),\label{eq2.3}\\
u_a(x)&=O(|x|^{-(N-2s)}) \quad\text{as }|x|\to\infty.\label{eq2.4}
\end{align}
In particular,
\[
u_a\in L^{\sstar}(\R^N\setminus B_R)
\]
for every sufficiently large $R$. If $x_0\in\Gamma$ is a nonremovable singular point of $u$, then $\mathcal I_a(x_0)$ is a nonremovable singular point of $u_a$.
\end{proposition}

\begin{proof}
The capacity and hyperplane assertions follow from Lemma~\ref{Lem2.2}. Since the inversion is a smooth diffeomorphism away from $a$, positivity, continuity, and local $H^s$-regularity are preserved away from $\Gamma_a$. It remains to prove that $u_a\in\cL_s(\R^N)$. Near $a$, use the change of variables $y=\mathcal I_a(x)$ and \eqref{eq2.1}:
\[
\int_{B_\rho(a)}|u_a(x)|\dd x = \int_{\R^N\setminus B_{1/\rho}(a)} |u(y)||y-a|^{-N-2s}\dd y<\infty.
\]
On compact subsets of $\R^N\setminus\{a\}$, local integrability is immediate. At infinity, $\mathcal I_a(x)$ remains in a compact neighborhood of the regular point $a$, and hence $u_a(x)\le C|x-a|^{-(N-2s)}$. Combining these estimates with the comparability of $|y-a|$ and $|y|$ at infinity shows that $u_a\in\cL_s(\R^N)$.

To prove \eqref{eq2.3}, use \eqref{eq2.1}, the identity
\[
(N-2s)\hstar=2N-\mu,
\]
and the change of variables $z=\mathcal I_a(y)$:
\begin{align*}
\cH[u_a](x) &=\int_{\R^N}\frac{|y-a|^{-(2N-\mu)}u^{\hstar}(\mathcal I_a(y))}{|x-y|^\mu}\dd y\\
&=|x-a|^{-\mu}\int_{\R^N}\frac{u^{\hstar}(z)}{|\mathcal I_a(x)-z|^\mu}\dd z =|x-a|^{-\mu}\cH[u](\mathcal I_a(x)).
\end{align*}
Thus the transformed potential is finite almost everywhere and the transformed source is locally integrable away from $\Gamma_a$.

Let $\phi\in C_c^\infty(\R^N\setminus\Gamma_a)$ and define
\[
\phi_a(y)=|y-a|^{-(N-2s)}\phi(\mathcal I_a(y)).
\]
Then $\phi_a\in C_c^\infty(\R^N\setminus\Gamma)$. Compact support follows because $\phi$ vanishes both near $a$ and outside a bounded set, while inversion exchanges neighborhoods of $a$ and infinity. We verify the Kelvin covariance directly; see also \cite[Proposition A.1]{RosOtonSerra2014}. Put $m=N-2s$, $w_a(z)=|z-a|^{-m}$, and $\xi=\mathcal I_a(y)$. The fractional fundamental solution satisfies $(-\Delta)^s w_a=0$ away from $a$. Define
\[
b(z)=\phi(z)-\phi(\xi)|\xi-a|^m w_a(z),\qquad b(\xi)=0.
\]
Substitution $t=\mathcal I_a(z)$ in the principal-value integral and \eqref{eq2.1} give
\begin{align*}
(-\Delta)^s\phi_a(y)&=C_{N,s}|y-a|^{-N-2s}\operatorname{PV}\int_{\R^N}\frac{-b(z)}{|\xi-z|^{N+2s}}\dd z\\
&=|y-a|^{-N-2s}(-\Delta)^s b(\xi).
\end{align*}
For completeness, the exclusion $|t-y|>\eps$ becomes $|z-\xi|>\eps|\xi-a||z-a|$. With $r=\eps|\xi-a|^2$, its symmetric difference from $|z-\xi|>r$ is contained in a shell of thickness $O(r^2)$. Since $b(z)=O(|z-\xi|)$ near $\xi$, changing the truncation produces an error $O(r^{2-2s})\to0$. This justifies the principal-value substitution for every $s\in(0,1)$. The identity $(-\Delta)^s b(\xi)=(-\Delta)^s\phi(\xi)$ therefore yields
\[
(-\Delta)^s\phi_a(y)=|y-a|^{-(N+2s)}((-\Delta)^s\phi)(\mathcal I_a(y)).
\]
Changing variables and using the weak equation for $u$, we obtain
\begin{align*}
\int_{\R^N}u_a(x)(-\Delta)^s\phi(x)\dd x &=\int_{\R^N}u(y)(-\Delta)^s\phi_a(y)\dd y\\
&=\int_{\R^N}\cH[u](y)u^{\hstar-1}(y)\phi_a(y)\dd y\\
&=\int_{\R^N}\cH[u_a](x)u_a^{\hstar-1}(x)\phi(x)\dd x.
\end{align*}
In the last step the exponent cancellation is
\[
\mu+(N-2s)\hstar-2N=0.
\]
For compactly supported smooth tests, the distributional pairing with $(-\Delta)^s$ equals the bilinear form used in Definition~\ref{Def1.1}. This proves \eqref{eq2.2}.

If $x_0\in\Gamma$ is nonremovable and $x_j\to x_0$ with $u(x_j)\to\infty$, then $\mathcal I_a(x_j)\to \mathcal I_a(x_0)$ and
\[
u_a(\mathcal I_a(x_j))=|x_j-a|^{N-2s}u(x_j)\to\infty.
\]
Hence $\mathcal I_a(x_0)$ is nonremovable. Finally, $\mathcal I_a(x)\to a$ as $|x|\to\infty$. Since $a$ is a regular point and $u$ is continuous there, \eqref{eq2.4} follows. Raising \eqref{eq2.4} to the power $\sstar$ gives $u_a^{\sstar}(x)=O(|x|^{-2N})$.
\end{proof}

\begin{remark}\label{Rem2.4}
The following example illustrates the decay obtained by Kelvin inversion at a regular point. Set $\gamma=(N-2s)/2$ and define $c_1,c_2>0$ by the fractional power and Riesz-potential identities
\[
(-\Delta)^s|x|^{-\gamma}=c_1|x|^{-\gamma-2s},\qquad
\int_{\R^N}\frac{|y|^{-(N-\mu/2)}}{|x-y|^\mu}\dd y=c_2|x|^{-\mu/2},\qquad x\ne0.
\]
Here $0<\gamma<N-2s$ gives $c_1>0$, and the second integral converges at $0$, $x$, and infinity because $0<\mu<N$. Since $\gamma\hstar=N-\mu/2$ and $\mu/2+\gamma(\hstar-1)=\gamma+2s$, the constant
\[
A_*=\left(\frac{c_1}{c_2}\right)^{1/(2\hstar-2)}
\]
makes
\[
U_0(x)=A_*|x|^{-\frac{N-2s}{2}}
\]
the homogeneous critical solution with singular set $\{0\}$; see also \cite[Theorem 1.3]{GuJevnikarYang2026a}. For $a\in\{x_1=0\}\setminus\{0\}$, its Kelvin transform is
\[
U_a(x)=|x-a|^{-(N-2s)}U_0(\mathcal I_a(x)) =A_*|a|^{-\frac{N-2s}{2}} \bigl(|x-a|\,|x-\mathcal I_a(0)|\bigr)^{-\frac{N-2s}{2}}.
\]
This function also solves the critical equation. Its singular set is
\[
\{a,\mathcal I_a(0)\}\subset\{x_1=0\}.
\]
The solution has the critical singular rates $|x-a|^{-(N-2s)/2}$ and $|x-\mathcal I_a(0)|^{-(N-2s)/2}$ at the two singular points and satisfies $U_a(x)=O(|x|^{-(N-2s)})$ at infinity. Thus $U_a\in L^{\sstar}(\R^N\setminus B_R)$ for every sufficiently large $R$, whereas $U_0$ does not belong to $L^{\sstar}$ near zero or at infinity.
\end{remark}

We next prove positivity and local coefficient bounds, and then obtain both critical tail estimates for every regular-center Kelvin transform. For $\delta>0$, write
\[
{O}_\delta(\Gamma)=\{x\in\R^N:\dist(x,\Gamma)<\delta\}.
\]

\begin{lemma}\label{Lem2.5}
Let $u$ be a positive weak singular solution of \eqref{eq1.1}. Then $u$ is strictly positive at every regular point and has a positive minimum on every nonempty compact subset of $\R^N\setminus\Gamma$. The potential $\cH[u]$ is finite at every regular point, and $\cH[u]$, the source $\cH[u]u^{\hstar-1}$, and the coefficient $\cA_u$ are locally bounded in $\R^N\setminus\Gamma$ for every $0<\mu<N$.
\end{lemma}

\begin{proof}
Choose $x_*\in\R^N\setminus\Gamma$ such that $\cH[u](x_*)<\infty$. For every $R>0$,
\[
\int_{B_R}u^{\hstar}(y)\dd y \le (R+|x_*|)^\mu\cH[u](x_*)<\infty.
\]
Fix $K\Subset\R^N\setminus\Gamma$ and choose $\rho>0$ such that the $2\rho$-neighborhood of $K$ is compactly contained in $\R^N\setminus\Gamma$. The function $u$ is uniformly bounded on that neighborhood by continuity. Choose $L$ sufficiently large that $B_\rho(x)\subset B_L$ for all $x\in K$. Split the Hartree integral into $B_\rho(x)$, $B_L\setminus B_\rho(x)$, and $\R^N\setminus B_L$. The first part is bounded because $\mu<N$ and $u$ is locally bounded. The second is bounded by the local $L^1$-integrability of $u^{\hstar}$. For the last part, if the ball is large enough, $|x-y|$ and $|x_*-y|$ are uniformly comparable, so it is controlled by $\cH[u](x_*)$. This proves local boundedness.

The source is locally bounded and nonnegative because $u$ is locally bounded and $\hstar-1>0$. Continuity and almost-everywhere positivity give $u\ge0$ at every regular point; set $u=0$ on the null set $\Gamma$.

Fix a regular point $z$ and a ball $B\Subset\R^N\setminus\Gamma$ containing $z$. Almost-everywhere positivity and continuity provide a ball $D\Subset\R^N\setminus(\Gamma\cup\overline B)$ and $d>0$ such that $u\ge d$ on $D$. Choose nonnegative functions $\varphi\in C_c^\infty(B)$ and $0\not\equiv\psi\in C_c^\infty(D)$ with $\varphi(z)=1$. For $x\in\overline B$,
\[
(-\Delta)^s\psi(x)
=-C_{N,s}\int_D\frac{\psi(y)}{|x-y|^{N+2s}}\dd y\le-c<0.
\]
The constant is uniform because the two supports are bounded and separated. The fractional Laplacian of $\varphi$ is bounded, as follows from its symmetric second differences and $s<1$. Choose $M>0$ so large that $h=\varphi+M\psi$ satisfies $(-\Delta)^sh\le0$ in $B$, and choose $\eps>0$ so small that $\eps M\|\psi\|_\infty\le d$.

Then $w=u-\eps h\ge0$ almost everywhere outside $B$ and $\cE(w,\phi)\ge0$ for every nonnegative $\phi\in C_c^\infty(B)$. Its negative part $w^-$ vanishes outside $\supp\varphi\Subset B$. If $0\le\chi\in C_c^\infty(B)$ equals one near $\supp\varphi$, then $w^-=(\chi w)^-\in H^s(\R^N)$ by the Lipschitz truncation inequality. The localization fact above and the local source bound justify testing with $w^-$, giving
\[
0\le\cE(w,w^-)
\le-\frac{C_{N,s}}2[w^-]_{H^s(\R^N)}^2.
\]
Here the second inequality follows by integrating
$(a-b)(a^--b^-)\le-|a^--b^-|^2$; the pairing is absolutely convergent. Hence $w^-=0$. Continuity at $z$ now gives $u(z)\ge\eps h(z)=\eps>0$.

Thus $u$ has a positive minimum on every regular compact set. The factor $u^{\hstar-2}$ is bounded on each such set regardless of the sign of $\hstar-2$. The local bound for $\cH[u]$ therefore proves the assertion for $\cA_u$.
\end{proof}

\Needspace{6\baselineskip}
\begin{proposition}\label{Prop2.6}
Let $u$ be a positive weak singular solution of \eqref{eq1.1}.
\begin{enumerate}
\item[(i)] If, for some $R_0>0$,
\begin{equation}\label{eq2.5}
u\in L^{\sstar}(\R^N\setminus B_{R_0}),
\qquad \cA_u\in L^{N/(2s)}(\R^N\setminus B_{R_0}),
\end{equation}
then, for every $\delta>0$,
\[
\begin{split}
u&\in L^{\sstar}(\R^N\setminus O_\delta(\Gamma)),\\
\cH[u]&\in L^{2N/\mu}(\R^N\setminus O_\delta(\Gamma)),\qquad
\cA_u\in L^{N/(2s)}(\R^N\setminus O_\delta(\Gamma)),
\end{split}
\]
and
\begin{equation}\label{eq2.6}
\|u\|_{L^{\sstar}(\R^N\setminus B_R)}
+\|\cA_u\|_{L^{N/(2s)}(\R^N\setminus B_R)}\longrightarrow0
\quad\text{as }R\to\infty.
\end{equation}
\item[(ii)] Without assuming \eqref{eq2.5} for $u$, let $a\in\{x_1=0\}\setminus\Gamma$. The Kelvin transform $u_a$ of Proposition~\ref{Prop2.3} satisfies \eqref{eq2.5}. More precisely,
\begin{equation}\label{eq2.AKelvin}
\cA_{u_a}(x)=|x-a|^{-4s}\cA_u(\mathcal I_a(x)),
\qquad x\in\R^N\setminus\Gamma_a,
\end{equation}
and, for every sufficiently large $L$,
\begin{equation}\label{eq2.tailidentity}
\begin{aligned}
\int_{\R^N\setminus B_L(a)}u_a^{\sstar}(x)\dd x
&=\int_{B_{1/L}(a)}u^{\sstar}(y)\dd y,\\
\int_{\R^N\setminus B_L(a)}\cA_{u_a}^{N/(2s)}(x)\dd x
&=\int_{B_{1/L}(a)}\cA_u^{N/(2s)}(y)\dd y.
\end{aligned}
\end{equation}
In particular, both integrals tend to zero as $L\to\infty$.
\end{enumerate}
\end{proposition}

\begin{proof}
For (i), Lemma~\ref{Lem2.5} bounds $u$ and $\cA_u$ on bounded sets whose closures are disjoint from $\Gamma$. Together with \eqref{eq2.5}, this proves their asserted integrability away from $O_\delta(\Gamma)$ and the tail limit \eqref{eq2.6}.

To estimate the potential, fix $\delta>0$ and write
\[
f_0=u^{\hstar}\one_{O_{\delta/2}(\Gamma)},\qquad
f_1=u^{\hstar}\one_{\R^N\setminus O_{\delta/2}(\Gamma)}.
\]
The first function has bounded support and belongs to $L^1$ by the local integrability proved in Lemma~\ref{Lem2.5}. Set
\[
r_0=\frac{2N}{2N-\mu}.
\]
Since $r_0\hstar=\sstar$, the first part of \eqref{eq2.5} and local boundedness give $f_1\in L^{r_0}(\R^N)$. The Hardy--Littlewood--Sobolev inequality yields
\[
|\cdot|^{-\mu}*f_1\in L^{2N/\mu}(\R^N).
\]
For $x\notin O_\delta(\Gamma)$, the potential generated by $f_0$ is bounded on bounded sets, by separation from its support, and is $O(|x|^{-\mu})$ at infinity. It therefore belongs to $L^{2N/\mu}(\R^N\setminus O_\delta(\Gamma))$, proving the assertion for $\cH[u]$.

For (ii), strict positivity from Lemma~\ref{Lem2.5} makes the coefficient well-defined at every regular point. The potential covariance \eqref{eq2.3} and
\[
\mu+(N-2s)(\hstar-2)=4s
\]
give \eqref{eq2.AKelvin} directly. Under the change of variables $y=\mathcal I_a(x)$, the Jacobian is $|y-a|^{-2N}$ and $|x-a|=|y-a|^{-1}$. The identities
\[
(N-2s)\sstar=2N,\qquad 4s\,\frac{N}{2s}=2N
\]
cancel the weights in the two integrals and prove \eqref{eq2.tailidentity}; the omitted center on the right has measure zero. For sufficiently large $L$, $B_{1/L}(a)$ lies in a fixed regular neighborhood of $a$. Both $u$ and $\cA_u$ are bounded there, so the right-hand sides are $O(L^{-N})$. In particular,
\[
\|u_a\|_{L^{\sstar}(\R^N\setminus B_L(a))}
=O(L^{-(N-2s)/2}),\qquad
\|\cA_{u_a}\|_{L^{N/(2s)}(\R^N\setminus B_L(a))}
=O(L^{-2s}).
\]
For $L>2|a|$, the inclusion $\R^N\setminus B_L\subset\R^N\setminus B_{L/2}(a)$ gives the same tail conclusions for balls centered at the origin. This proves \eqref{eq2.5} for $u_a$, independently of any such hypothesis on $u$.
\end{proof}

\begin{lemma}\label{Lem2.7}
Let $E\subset\R^N$ be measurable, let $u\in L^{\sstar}(E)$ be nonnegative, and let $Z\in L^{\sstar}(\R^N)$ be nonnegative and vanish almost everywhere outside $E$. Then
\begin{equation}\label{eq2.7}
\iint_{E\times E} \frac{u^{\hstar-1}(x)Z(x)u^{\hstar-1}(y)Z(y)}{|x-y|^\mu}\dd x\dd y \le C\|u\|_{L^{\sstar}(E)}^{2\hstar-2}\|Z\|_{L^{\sstar}(\R^N)}^2.
\end{equation}
If, in addition, $\cA_u\in L^{N/(2s)}(E)$, then
\begin{equation}\label{eq2.8}
\int_E\cA_u Z^2\dd x \le \|\cA_u\|_{L^{N/(2s)}(E)}\|Z\|_{L^{\sstar}(\R^N)}^2.
\end{equation}
\end{lemma}

\begin{proof}
Let $r_0=2N/(2N-\mu)$. Then $1<r_0<2$, $2/r_0+\mu/N=2$, and $1/r_0=\hstar/\sstar$. Since $\hstar>1$, the H\"older exponents $\hstar/(\hstar-1)$ and $\hstar$ are both greater than one, including when $\hstar<2$. The Hardy--Littlewood--Sobolev and H\"older inequalities give
\begin{align*}
\iint_{E\times E} \frac{u^{\hstar-1}(x)Z(x)u^{\hstar-1}(y)Z(y)}{|x-y|^\mu}\dd x\dd y &\le C\|u^{\hstar-1}Z\|_{L^{r_0}(E)}^2\\
&\le C\|u\|_{L^{\sstar}(E)}^{2\hstar-2} \|Z\|_{L^{\sstar}(\R^N)}^2,
\end{align*}
which proves \eqref{eq2.7}. Since $N/(2s)$ and $\sstar/2$ are conjugate exponents, H\"older's inequality also gives \eqref{eq2.8}.
\end{proof}

\begin{lemma}\label{Lem2.8}
For every $q>0$ and $0\le r<t$,
\begin{equation}\label{eq2.power}
t^q-r^q\le c_qt^{q-1}(t-r),\qquad c_q=\max\{1,q\}.
\end{equation}
\end{lemma}

\begin{proof}
If $q\ge1$, integrate $q\tau^{q-1}\le qt^{q-1}$ over $[r,t]$.
If $0<q\le1$, set $z=r/t\in[0,1)$ and use $z^q\ge z$ to obtain
$t^q-r^q=t^q(1-z^q)\le t^q(1-z)=t^{q-1}(t-r)$.
\end{proof}

\section{Truncated Reflected Comparison}\label{Sec3}

Throughout this section, $\Gamma\subset\{x_1=0\}$ is compact with zero $(s,2)$-capacity, and $u$ is a positive weak singular solution. 
For $\lambda\in\R$, set
\[
T_\lambda=\{x_1=\lambda\}, \qquad \Sigma_\lambda=\{x_1<\lambda\}, \qquad x^\lambda=(2\lambda-x_1,x'),
\]
\[
u_\lambda(x)=u(x^\lambda), \qquad \Gamma_\lambda=\{x^\lambda:x\in\Gamma\}, \qquad w_\lambda=u-u_\lambda.
\]
For $\lambda<0$, define
\[
W_\lambda=(w_\lambda)^+\one_{\Sigma_\lambda}, \qquad E_\lambda=\{x\in\Sigma_\lambda\setminus\Gamma_\lambda:u(x)>u_\lambda(x)\}.
\]
Since $\Gamma\subset\{x_1=0\}$, the half-space $\Sigma_\lambda$ has positive distance from $\Gamma$, while $\Gamma_\lambda\subset\Sigma_\lambda$.

\begin{lemma}\label{Lem3.1}
Let $0\le\xi\le1$ be a smooth compactly supported function that vanishes in a neighborhood of $\Gamma_\lambda$. Then
\[
W_\lambda\xi\in H^s(\R^N), \qquad W_\lambda\xi^2\in H_0^s(\Omega)
\]
for some bounded open set
\[
\Omega\Subset\R^N\setminus(\Gamma\cup\Gamma_\lambda).
\]
Moreover,
\begin{equation}\label{eq3.1}
\cE(w_\lambda,W_\lambda\xi^2) = \int_{E_\lambda} \left(\cH[u]u^{\hstar-1}-\cH[u_\lambda]u_\lambda^{\hstar-1}\right) W_\lambda\xi^2\dd x.
\end{equation}
\end{lemma}

\begin{proof}
The reflected function $u_\lambda$ belongs to $\cL_s(\R^N)\cap H^s_{\loc}(\R^N\setminus\Gamma_\lambda)$ and solves the reflected equation. Indeed, reflection preserves Lebesgue measure and Euclidean distance, while
\[
1+|x^\lambda|^{N+2s}\asymp_\lambda 1+|x|^{N+2s}.
\]
Moreover,
\[
\cH[u_\lambda](x)=\cH[u](x^\lambda).
\]
If $\psi\in C_c^\infty(\R^N\setminus\Gamma_\lambda)$ and $\psi^\lambda(x)=\psi(x^\lambda)$, then $(-\Delta)^s\psi^\lambda(x)=((-\Delta)^s\psi)(x^\lambda)$. Testing the equation of $u$ with $\psi^\lambda$ and changing variables proves the weak equation for $u_\lambda$.

Let
\[
K=\supp\xi\cap\overline{\Sigma_\lambda}.
\]
Because $\xi$ vanishes near $\Gamma_\lambda$ and $\Sigma_\lambda$ has positive distance from $\Gamma$, the compact set $K\cup K^\lambda$ is disjoint from $\Gamma\cup\Gamma_\lambda$. Since the complement of $\Gamma\cup\Gamma_\lambda$ is invariant under reflection across $T_\lambda$, we may choose
\[
\chi\in C_c^\infty(\R^N\setminus(\Gamma\cup\Gamma_\lambda))
\]
such that $\chi(x^\lambda)=\chi(x)$ and $\chi=1$ near $K\cup K^\lambda$. Then
\[
g=\chi w_\lambda\in H^s(\R^N), \qquad g(x^\lambda)=-g(x).
\]
We record the antisymmetric truncation estimate. If $g\in H^s(\R^N)$ is antisymmetric and $G=g^+\one_{\Sigma_\lambda}$, set
\begin{align*}
I_1&=\iint_{\Sigma_\lambda\times\Sigma_\lambda}\frac{|g(x)-g(y)|^2}{|x-y|^{N+2s}}\dd x\dd y,\\
I_2&=\iint_{\Sigma_\lambda\times\Sigma_\lambda}\frac{|g(x)+g(y)|^2}{|x-y^\lambda|^{N+2s}}\dd x\dd y.
\end{align*}
Reflection and antisymmetry give
\[
[g]_{H^s(\R^N)}^2=2I_1+2I_2.
\]
The contribution to $[G]_{H^s}^2$ from $\Sigma_\lambda\times\Sigma_\lambda$ is at most $I_1$ by the Lipschitz continuity of the positive part map. For the cross interaction, reflect the exterior point into $y\in\Sigma_\lambda$. When $g(x)>0$,
\[
g^2(x)\le
\begin{cases}
|g(x)+g(y)|^2,&g(y)\ge0,\\
|g(x)-g(y)|^2,&g(y)<0.
\end{cases}
\]
In the second case use $|x-y|<|x-y^\lambda|$. The cross interaction, including both orders of the variables, is therefore at most $2I_1+2I_2$. Tonelli's theorem yields
\[
[G]_{H^s(\R^N)}^2\le3I_1+2I_2\le\frac32[g]_{H^s(\R^N)}^2.
\]
Since $|G|\le|g|$, this also proves $G\in H^s(\R^N)$. The same truncation is used in the weak antisymmetric maximum principle of \cite{FallJarohs2015}. On $K$, one has $G=W_\lambda$. Hence $W_\lambda\xi=G\xi\in H^s(\R^N)$, and multiplication by $\xi$ gives $W_\lambda\xi^2\in H^s(\R^N)$.

The support of $W_\lambda\xi^2$ is compact and disjoint from $\Gamma\cup\Gamma_\lambda$. It is therefore compactly contained in a bounded smooth open set $\Omega\Subset\R^N\setminus(\Gamma\cup\Gamma_\lambda)$, and the function belongs to $H_0^s(\Omega)$. The bilinear pairings with $u$ and $u_\lambda$ are continuous on $H_0^s(\Omega)$. Indeed, on a slightly larger bounded set they are controlled by local $H^s$-regularity, while the complementary interactions are controlled by the $\cL_s$-norm. By Lemma~\ref{Lem2.5}, both Hartree potentials are locally bounded. Since $\hstar>1$ and $u,u_\lambda$ are locally bounded away from their singular sets, both sources are locally bounded there as well. Approximation of $W_\lambda\xi^2$ by functions in $C_c^\infty(\Omega)$ and subtraction of the two weak equations prove \eqref{eq3.1}. The same local and tail estimates show that all bilinear pairings in this identity are absolutely convergent.
\end{proof}

Fix $\theta\in C_c^\infty(B_2)$ such that $0\le\theta\le1$ and $\theta=1$ on $B_1$, and put $\theta_R(x)=\theta(x/R)$. Choose decreasing neighborhoods $O_j$ of $\Gamma_\lambda$ such that
\[
\overline{O_{j+1}}\subset O_j, \qquad \overline{O_j}\Subset\Sigma_\lambda, \qquad \bigcap_{j=1}^\infty\overline{O_j}=\Gamma_\lambda.
\]
By Lemma~\ref{Lem2.1}, there are $\eta_j\in C^\infty(\R^N)$ such that
\[
0\le\eta_j\le1, \qquad \eta_j=0\text{ near }\Gamma_\lambda, \qquad \eta_j=1\text{ on }\R^N\setminus O_j, \qquad \|1-\eta_j\|_{H^s}<\frac1j.
\]
Set
\[
\xi_{j,R}=\eta_j\theta_R.
\]

\begin{lemma}\label{Lem3.2}
Assume the first condition in \eqref{eq2.5}. For the preceding cutoffs, define
\[
\mathcal R_{j,R}= \iint_{\Sigma_\lambda\times\Sigma_\lambda} \frac{W_\lambda(x)W_\lambda(y)|\xi_{j,R}(x)-\xi_{j,R}(y)|^2}{|x-y|^{N+2s}}\dd x\dd y.
\]
Then $\mathcal R_{j,R}<\infty$ and
\begin{equation}\label{eq3.2}
\lim_{R\to\infty}\limsup_{j\to\infty}\mathcal R_{j,R}=0.
\end{equation}
Moreover, $\eta_j\to1$ almost everywhere in $\R^N$.
\end{lemma}

\begin{proof}
The function $W_\lambda$ is bounded on bounded subsets of $\overline{\Sigma_\lambda}$, since $0\le W_\lambda\le u$ there and these sets have positive distance from $\Gamma$. Together with \eqref{eq2.5} and $W_\lambda=0$ outside $\Sigma_\lambda$, this gives $W_\lambda\in L^{\sstar}(\R^N)$.

Put $\zeta_j=1-\eta_j$. These functions are supported in a fixed bounded open set $O$ with $\overline O\Subset\Sigma_\lambda$, and $\|\zeta_j\|_{H^s}\to0$. Define
\[
J_j(x)=\int_{\R^N}\frac{|\zeta_j(x)-\zeta_j(y)|^2}{|x-y|^{N+2s}}\dd y.
\]
Using $2ab\le a^2+b^2$ and symmetry, we obtain
\[
\iint_{\Sigma_\lambda\times\Sigma_\lambda} \frac{W_\lambda(x)W_\lambda(y)|\eta_j(x)-\eta_j(y)|^2}{|x-y|^{N+2s}}\dd x\dd y \le\int_{\Sigma_\lambda}W_\lambda^2J_j\dd x.
\]
Choose a bounded open set $O'$ with $\overline O\Subset O'\Subset\Sigma_\lambda$. Since $W_\lambda$ is bounded on $O'$,
\[
\int_{O'}W_\lambda^2J_j\dd x \le\|W_\lambda\|_{L^\infty(O')}^2[\zeta_j]_{H^s(\R^N)}^2\longrightarrow0.
\]
For $x\notin O'$, separation from $\supp\zeta_j$ gives
\[
J_j(x)\le C\|\zeta_j\|_{L^2}^2(1+|x|)^{-N-2s}.
\]
The integral of $W_\lambda^2(1+|x|)^{-N-2s}$ is finite by H\"older's inequality and $W_\lambda\in L^{\sstar}(\R^N)$. Thus the capacitary contribution tends to zero.

For the far field cutoff, set
\[
J_R(x)=\int_{\R^N}\frac{|\theta_R(x)-\theta_R(y)|^2}{|x-y|^{N+2s}}\dd y =R^{-2s}J_1(x/R).
\]
The function $J_1$ is bounded and satisfies $J_1(x)\le C|x|^{-N-2s}$ for $|x|>4$. Hence
\[
\sup_{R>1}\|J_R\|_{L^{N/(2s)}(\R^N)}<\infty,
\]
and $J_R\to0$ uniformly on bounded sets. For each $L>0$, H\"older's inequality gives
\[
\int_{\Sigma_\lambda\setminus B_L}W_\lambda^2J_R\dd x \le C\|W_\lambda\|_{L^{\sstar}(\R^N\setminus B_L)}^2.
\]
This is uniformly small for large $L$, while the integral over $B_L$ tends to zero by local $L^2$-integrability of $W_\lambda$. Therefore
\[
\int_{\Sigma_\lambda}W_\lambda^2J_R\dd x\longrightarrow0 \quad\text{as }R\to\infty.
\]
The inequality
\[
|\eta_j(x)\theta_R(x)-\eta_j(y)\theta_R(y)|^2 \le2|\eta_j(x)-\eta_j(y)|^2+2|\theta_R(x)-\theta_R(y)|^2
\]
now proves finiteness of $\mathcal R_{j,R}$ and \eqref{eq3.2}. Finally, the nested neighborhoods $O_j$ eventually exclude every point outside $\Gamma_\lambda$, so $\eta_j=1$ there for all sufficiently large $j$. Thus $\eta_j\to1$ outside $\Gamma_\lambda$, and hence almost everywhere because $|\Gamma_\lambda|=0$.
\end{proof}

\begin{lemma}\label{Lem3.3}
Let $w$ be measurable, finite almost everywhere, and antisymmetric with respect to $T_\lambda$, let $W=w^+\one_{\Sigma_\lambda}$, and let $\xi$ be measurable with $0\le\xi\le1$. Assume that $W\xi\in H^s(\R^N)$, that
\[
\iint_{\Sigma_\lambda\times\Sigma_\lambda} \frac{W(x)W(y)|\xi(x)-\xi(y)|^2}{|x-y|^{N+2s}}\dd x\dd y<\infty,
\]
and that the pairing $\cE(w,W\xi^2)$ is absolutely convergent. Then
\begin{equation}\label{eq3.3}
\cE(w,W\xi^2) \ge \frac{C_{N,s}}2 \left([W\xi]_{H^s(\R^N)}^2 - \iint_{\Sigma_\lambda\times\Sigma_\lambda} \frac{W(x)W(y)|\xi(x)-\xi(y)|^2}{|x-y|^{N+2s}}\dd x\dd y\right).
\end{equation}
\end{lemma}

\begin{proof}
Write $\Sigma=\Sigma_\lambda$, $P=w^-\one_\Sigma$, and $A=W\xi^2$. Since $A=0$ on $\R^N\setminus\Sigma$, absolute convergence permits us to split the full bilinear form and reflect the second variable in the cross interaction. We obtain
\begin{align*}
\frac{2}{C_{N,s}}\cE(w,A) ={}&\iint_{\Sigma\times\Sigma} \frac{(w(x)-w(y))(A(x)-A(y))}{|x-y|^{N+2s}}\dd x\dd y\\
&+\iint_{\Sigma\times\Sigma} \frac{(w(x)+w(y))(A(x)+A(y))}{|x-y^\lambda|^{N+2s}}\dd x\dd y.
\end{align*}
To justify the following separation into nonnegative terms, apply the algebra first on the increasing sets
\[
\begin{aligned}
D_n=\bigl\{(x,y)\in(\Sigma\cap B_n)^2:\ & \min\{\dist(x,T_\lambda),\dist(y,T_\lambda)\}>\tfrac1n,\\
&|x-y|>\tfrac1n,\quad |w(x)|+|w(y)|\le n\bigr\}.
\end{aligned}
\]
All integrals on $D_n$ are finite. The original two integrands converge in $L^1$ along this exhaustion, while every term retained below is nonnegative except for the last integral in \eqref{eq3.3}, which is finite by assumption. Thus monotone convergence and Fatou's lemma allow passage to the full half-space.

On $\Sigma$, one has $w=W-P$ and $WP=0$. The terms containing both $P$ and $W$ combine into
\[
\bigl(P(x)W(y)\xi^2(y)+P(y)W(x)\xi^2(x)\bigr) \left( \frac1{|x-y|^{N+2s}}- \frac1{|x-y^\lambda|^{N+2s}} \right),
\]
which is nonnegative because $|x-y|<|x-y^\lambda|$ for $x,y\in\Sigma$. For the remaining same half-space term, set $Z=W\xi$ and use
\[
(W(x)-W(y))(W(x)\xi^2(x)-W(y)\xi^2(y)) =(Z(x)-Z(y))^2-W(x)W(y)(\xi(x)-\xi(y))^2.
\]
The $W$-part of the reflected cross interaction is bounded below by
\[
\frac{Z^2(x)+Z^2(y)}{|x-y^\lambda|^{N+2s}}.
\]
Consequently,
\begin{align*}
\frac{2}{C_{N,s}}\cE(w,A) \ge{}&\iint_{\Sigma\times\Sigma} \frac{(Z(x)-Z(y))^2}{|x-y|^{N+2s}}\dd x\dd y\\
&+\iint_{\Sigma\times\Sigma} \frac{Z^2(x)+Z^2(y)}{|x-y^\lambda|^{N+2s}}\dd x\dd y\\
&-\iint_{\Sigma\times\Sigma} \frac{W(x)W(y)|\xi(x)-\xi(y)|^2}{|x-y|^{N+2s}}\dd x\dd y.
\end{align*}
After the change of variables $z=y^\lambda$, the first two integrals on the right are exactly $[Z]_{H^s(\R^N)}^2$. This proves \eqref{eq3.3}.
\end{proof}

Let
\[
V=\cH[u], \qquad V_\lambda(x)=V(x^\lambda)=\cH[u_\lambda](x).
\]

\begin{lemma}\label{Lem3.4}
For every $x\in\Sigma_\lambda\setminus\Gamma_\lambda$,
\begin{equation}\label{eq3.4}
V(x)-V_\lambda(x) = \int_{\Sigma_\lambda} \left(u^{\hstar}(y)-u_\lambda^{\hstar}(y)\right) \left(\frac1{|x-y|^\mu}-\frac1{|x^\lambda-y|^\mu}\right)\dd y.
\end{equation}
Moreover, for distinct $x,y\in\Sigma_\lambda$,
\[
\frac1{|x-y|^\mu}- \frac1{|x^\lambda-y|^\mu}\ge0,
\]
and for every $x\in\Sigma_\lambda\setminus\Gamma_\lambda$,
\begin{equation}\label{eq3.5}
(V(x)-V_\lambda(x))^+ \le \hstar\int_{E_\lambda} \frac{u^{\hstar-1}(y)W_\lambda(y)}{|x-y|^\mu}\dd y.
\end{equation}
\end{lemma}

\begin{proof}
The sets $\Gamma$ and $\Gamma_\lambda$ have measure zero. To obtain \eqref{eq3.4}, split the integrals defining $V(x)$ and $V(x^\lambda)$ into the two half-spaces and reflect the complementary half-space. The kernel sign follows from
\[
|x^\lambda-y|^2-|x-y|^2 =4(\lambda-x_1)(\lambda-y_1)\ge0.
\]
The integrand in \eqref{eq3.4} is absolutely integrable because the kernel difference is bounded by $|x-y|^{-\mu}$ and the resulting integrals are bounded by the finite potentials $V(x)$ and $V_\lambda(x)$. On $\Sigma_\lambda\setminus E_\lambda$, the density difference is nonpositive. On $E_\lambda$,
\[
u^{\hstar}-u_\lambda^{\hstar} \le\hstar u^{\hstar-1}W_\lambda.
\]
This proves \eqref{eq3.5}.
\end{proof}

For measurable $E\subset\R^N\setminus\Gamma$, define
\[
\mathfrak M_u(E)=\|u\|_{L^{\sstar}(E)}^{2\hstar-2}+\|\cA_u\|_{L^{N/(2s)}(E)}.
\]

\begin{proposition}\label{Prop3.5}
Assume the two tail conditions \eqref{eq2.5}. There is a constant $\kappa_*>0$, depending only on $N,s,\mu$, such that, for every $\lambda<0$,
\[
\mathfrak M_u(E_\lambda)<\kappa_* \quad\Longrightarrow\quad u(x)\le u_\lambda(x) \quad\text{for every }x\in\Sigma_\lambda\setminus\Gamma_\lambda.
\]
The same constant applies when the planes move from the right.
\end{proposition}

\begin{proof}
By Proposition~\ref{Prop2.6}, $W_\lambda\in L^{\sstar}(\R^N)$ and $\mathfrak M_u(E_\lambda)<\infty$. Set $Z_{j,R}=W_\lambda\xi_{j,R}$. Lemma~\ref{Lem3.1} gives $Z_{j,R}\in H^s(\R^N)$ and
\[
\cE(w_\lambda,W_\lambda\xi_{j,R}^2) =\int_{E_\lambda} (Vu^{\hstar-1}-V_\lambda u_\lambda^{\hstar-1})W_\lambda\xi_{j,R}^2\dd x.
\]
On $E_\lambda$, Lemma~\ref{Lem2.8} with $q=\hstar-1>0$ gives
\[
\begin{aligned}
Vu^{\hstar-1}-V_\lambda u_\lambda^{\hstar-1}
&=V(u^{\hstar-1}-u_\lambda^{\hstar-1})
+(V-V_\lambda)u_\lambda^{\hstar-1}\\
&\le c_{\hstar-1}\cA_uW_\lambda
+(V-V_\lambda)^+u^{\hstar-1}.
\end{aligned}
\]
The coefficient is evaluated at the larger value $u$, also when $\hstar<2$.
Since $0\le\xi_{j,R}\le1$, Lemma~\ref{Lem3.4} yields
\begin{align*}
\cE(w_\lambda,W_\lambda\xi_{j,R}^2) \le{}&C\int_{E_\lambda}\cA_uW_\lambda^2\dd x\\
&+C\iint_{E_\lambda\times E_\lambda} \frac{u^{\hstar-1}(x)W_\lambda(x)u^{\hstar-1}(y)W_\lambda(y)}{|x-y|^\mu}\dd x\dd y.
\end{align*}
Both terms are finite by Lemma~\ref{Lem2.7}, applied with $Z=W_\lambda$. It follows that
\begin{equation}\label{eq3.6}
\cE(w_\lambda,W_\lambda\xi_{j,R}^2) \le C\mathfrak M_u(E_\lambda)\|W_\lambda\|_{L^{\sstar}(\R^N)}^2.
\end{equation}
Since $\mathcal R_{j,R}<\infty$ by Lemma~\ref{Lem3.2}, we may apply Lemma~\ref{Lem3.3}:
\begin{equation}\label{eq3.7}
\cE(w_\lambda,W_\lambda\xi_{j,R}^2) \ge\frac{C_{N,s}}2\left([Z_{j,R}]_{H^s(\R^N)}^2-\mathcal R_{j,R}\right).
\end{equation}
Combining \eqref{eq3.6} and \eqref{eq3.7}, we obtain
\[
[Z_{j,R}]_{H^s(\R^N)}^2 \le C\mathfrak M_u(E_\lambda)\|W_\lambda\|_{L^{\sstar}(\R^N)}^2+\mathcal R_{j,R}.
\]
Since $Z_{j,R}\in H^s(\R^N)$ has compact support and $N>2s$, the fractional Sobolev inequality \cite[Theorem 6.5]{DiNezzaPalatucciValdinoci2012} applies to $Z_{j,R}$ and yields
\[
\|Z_{j,R}\|_{L^{\sstar}(\R^N)}^2\le C\mathfrak M_u(E_\lambda)\|W_\lambda\|_{L^{\sstar}(\R^N)}^2+C\mathcal R_{j,R}.
\]
By \eqref{eq3.2}, choose $R_n\ge n$ with $\limsup_{j\to\infty}\mathcal R_{j,R_n}<1/(2n)$, and then choose $j_n\ge n$ with $\mathcal R_{j_n,R_n}<1/n$. The shrinking neighborhoods $O_j$ and the identity $\theta_R=1$ on $B_R$ imply $Z_{j_n,R_n}\to W_\lambda$ almost everywhere. Fatou's lemma in $L^{\sstar}$ now gives
\[
\|W_\lambda\|_{L^{\sstar}(\R^N)}^2\le\liminf_{n\to\infty}\|Z_{j_n,R_n}\|_{L^{\sstar}(\R^N)}^2\le C\mathfrak M_u(E_\lambda)\|W_\lambda\|_{L^{\sstar}(\R^N)}^2.
\]
Choose $\kappa_*>0$ such that $C\kappa_*<1$. If $\mathfrak M_u(E_\lambda)<\kappa_*$, the finiteness of $\|W_\lambda\|_{L^{\sstar}}$ implies $W_\lambda=0$ almost everywhere. Since both functions are continuous in $\Sigma_\lambda\setminus\Gamma_\lambda$, the comparison holds at every point there. Reflecting the argument across $\{x_1=0\}$ proves the result for planes moving from the right, with the same constant.
\end{proof}

\section{Moving Planes for the Kelvin Transform}\label{Sec4}

For a bounded open set $\Omega\subset\R^N$, set
\[
\cD^s_{\loc}(\Omega)=\left\{ v\in\cL_s(\R^N): v\in H^s_{\loc}(O)\text{ for some open }O\text{ with }\Omega\Subset O \right\}.
\]
For $v\in\cD^s_{\loc}(\Omega)$ and $\phi\in H_0^s(\Omega)$, the pairing $\cE(v,\phi)$ is absolutely convergent and depends continuously on $\phi$ in $H_0^s(\Omega)$. To see this, choose $\chi\in C_c^\infty(O)$ equal to one in a neighborhood of $\overline\Omega$. Then $\chi v\in H^s(\R^N)$, so its pairing with $\phi$ is controlled by the Cauchy--Schwarz inequality. The remaining interaction involves $(1-\chi)v$, whose support has positive distance from $\Omega$, and is controlled by the $\cL_s$ condition. Thus $\cD^s_{\loc}(\Omega)$ is contained in the class $\mathcal D^s(\Omega)$ defined by finite pairings in \cite[Section 2]{FallJarohs2015}.

\begin{lemma}\label{Lem4.1}
Let $H$ be a half-space, let $B\Subset H$ be a ball, and let $v\in\cD^s_{\loc}(B)$. Let $Q$ denote reflection across $\partial H$. Assume that $v(Qx)=-v(x)$ for every $x\in\R^N$, that $v\ge0$ almost everywhere in $H$, and
\[
\cE(v,\phi)\ge0
\]
for every nonnegative $\phi\in H_0^s(B)$. Then either $v=0$ almost everywhere in $\R^N$, or $v>0$ almost everywhere in $B$. In the nonzero alternative, any continuous representative of $v$ in $B$ is strictly positive at every point of $B$.
\end{lemma}

\begin{proof}
We apply \cite[Corollary 3.4]{FallJarohs2015}, with $\Omega=B$ and $c\equiv0$. Set $\mathcal N=\{x\in H:v(x)<0\}$. Since $\mathcal N$ is null, redefining $v$ to be zero on $\mathcal N\cup Q(\mathcal N)$ leaves all pairings unchanged, preserves exact antisymmetry, and makes $v\ge0$ throughout $H$. The set $B$ is open and bounded, with $B\Subset H$. By the inclusion established above, $v\in\cD^s_{\loc}(B)$ belongs to the class $\mathcal D^s(B)$ used in that reference: the pairing $\cE(v,\phi)$ is finite for every $\phi\in H_0^s(B)$, and the assumed variational inequality says that
\[
(-\Delta)^s v\ge0 \quad\text{weakly in }B.
\]
Moreover, $v\ge0$ throughout $H$, in particular on $H\setminus B$, and $v(Qx)=-v(x)$ for every $x\in\R^N$. It follows from \cite[Corollary 3.4]{FallJarohs2015} that either $v=0$ almost everywhere in $H$ or $v>0$ almost everywhere in $B$. In the first case, antisymmetry gives $v=0$ almost everywhere in $\R^N$.

Assume now that the second alternative holds and let $\widetilde v\in C(B)$ be a continuous representative of $v$. Fix $x_0\in B$. Choose a ball $B'=B_r(x_0)$ with $\overline{B'}\subset B$, small enough that $B\setminus\overline{B'}$ contains a nonempty open set $D$ with $\overline D\subset B\setminus\overline{B'}$. Since $v>0$ almost everywhere in $D$, there exists $m\in\mathbb N$ such that
\[
K=\{x\in\overline D:\widetilde v(x)\ge m^{-1}\}
\]
has positive measure. Since $v=\widetilde v$ almost everywhere in $B$,
\[
K\Subset H\setminus\overline{B'}, \qquad |K|>0, \qquad \operatorname*{ess\,inf}_{K}v\ge m^{-1}>0.
\]
We apply \cite[Proposition 3.3]{FallJarohs2015} to $v$ with $\Omega=B$, the ball $B'$, this set $K$, and $c\equiv0$; indeed, $\|c\|_{L^\infty(B)}=0<\lambda_1(B')$. Here $\lambda_1(B')>0$ is the first Dirichlet eigenvalue of $(-\Delta)^s$ on $B'$. It follows that, for some $d>0$,
\[
v(x)\ge d\,\dist(x,\R^N\setminus B')^s \quad\text{for almost every }x\in B'.
\]
Since $v=\widetilde v$ almost everywhere in $B'$, continuity gives $\widetilde v(x)\ge d\,\dist(x,\R^N\setminus B')^s$ at every point of $B'$. In particular, $\widetilde v(x_0)\ge d r^s>0$. Since $x_0$ was arbitrary, $\widetilde v>0$ at every point of $B$.
\end{proof}

\begin{proposition}\label{Prop4.2}
Let $\lambda<0$ and suppose that
\[
u\le u_\lambda \quad\text{in }\Sigma_\lambda\setminus\Gamma_\lambda.
\]
Then either $u=u_\lambda$ almost everywhere in $\R^N$, or
\[
u(x)<u_\lambda(x) \quad\text{for every }x\in\Sigma_\lambda\setminus\Gamma_\lambda.
\]
The corresponding alternative holds for planes moving from the right.
\end{proposition}

\begin{proof}
Set $v_\lambda=u_\lambda-u$ on $\R^N\setminus(\Gamma\cup\Gamma_\lambda)$ and set $v_\lambda=0$ on $\Gamma\cup\Gamma_\lambda$. This null set is invariant under reflection across $T_\lambda$, so $v_\lambda(x^\lambda)=-v_\lambda(x)$ for every $x\in\R^N$ and $v_\lambda\ge0$ throughout $\Sigma_\lambda$. The choice of values on the null set does not affect any weak equation or energy pairing. Let $B\Subset\Sigma_\lambda\setminus\Gamma_\lambda$ be a ball, and choose a ball $\widehat B$ with $\overline B\subset\widehat B\Subset\Sigma_\lambda\setminus\Gamma_\lambda$. Local $H^s$-regularity and the $\cL_s$ condition imply $v_\lambda\in\cD^s_{\loc}(\widehat B)$, and hence $v_\lambda\in\cD^s_{\loc}(B)$.

Splitting the finite Hartree integrals as in \eqref{eq3.4}, the assumed comparison gives $V\le V_\lambda$ in $\Sigma_\lambda\setminus\Gamma_\lambda$. Since $\hstar>1$, it follows that
\[
V_\lambda u_\lambda^{\hstar-1}-Vu^{\hstar-1}\ge0.
\]
For nonnegative $\phi\in H_0^s(B)$, mollification of its zero extension gives an approximation in $H^s(\R^N)$ by nonnegative functions in $C_c^\infty(\widehat B)$. The continuity of the local energy pairing established above and the local boundedness of both sources from Lemma~\ref{Lem2.5} permit passage to the limit in the two weak equations. Hence
\[
\cE(v_\lambda,\phi) = \int_B (V_\lambda u_\lambda^{\hstar-1}-Vu^{\hstar-1})\phi\dd x \ge0.
\]
By Lemma~\ref{Lem4.1}, either $v_\lambda=0$ almost everywhere in $\R^N$, or $v_\lambda>0$ at every point of $B$. Every point of $\Sigma_\lambda\setminus\Gamma_\lambda$ lies in such a ball, so the proposition follows.
\end{proof}

Lemma~\ref{Lem4.1} and Proposition~\ref{Prop4.2} require neither of the critical tail conditions in \eqref{eq2.5}. The following obstruction to equality is also independent of them.

\begin{lemma}\label{Lem4.3}
If $x_0\in\Gamma$ is a nonremovable singular point, then, for every $\lambda\ne0$, the identity $u=u_\lambda$ cannot hold almost everywhere in $\R^N$.
\end{lemma}

\begin{proof}
Choose $x_j\to x_0$ in $\R^N\setminus\Gamma$ such that $u(x_j)\to\infty$. If $\lambda\ne0$, then $x_j^\lambda\to x_0^\lambda$, where $x_0^\lambda\notin\Gamma$ because its first coordinate is $2\lambda$. Since $\dist(\Gamma,\Gamma_\lambda)=2|\lambda|$, both $x_j\notin\Gamma_\lambda$ and $x_j^\lambda\notin\Gamma$ for all large $j$. If $u=u_\lambda$ almost everywhere, continuity away from the two singular sets makes the identity pointwise near $x_j^\lambda$, and
\[
u(x_j^\lambda)=u_\lambda(x_j^\lambda)=u(x_j)\to\infty.
\]
This contradicts local boundedness at the regular point $x_0^\lambda$.
\end{proof}

For the remainder of this section, assume the two tail conditions \eqref{eq2.5} and the existence of a nonremovable singular point on $\Gamma$. For $\lambda<0$, the set $\Sigma_\lambda$ has positive distance from $\Gamma$, and Proposition~\ref{Prop2.6} shows that $\mathfrak M_u(E)$ is finite for every measurable $E\subset\Sigma_\lambda$.

Since $\Sigma_\lambda\subset\R^N\setminus B_{|\lambda|}$ for $\lambda<0$, \eqref{eq2.6} gives
\[
\mathfrak M_u(\Sigma_\lambda)\to0 \quad\text{as }\lambda\to-\infty.
\]
Since $E_\lambda\subset\Sigma_\lambda$, Proposition~\ref{Prop3.5} gives
\[
u\le u_\lambda \quad\text{in }\Sigma_\lambda\setminus\Gamma_\lambda
\]
for every sufficiently negative $\lambda$.

Define
\[
\Lambda_-= \left\{ \lambda<0: u\le u_t\text{ in }\Sigma_t\setminus\Gamma_t \text{ for every }t\le\lambda \right\}, \qquad \lambda_0=\sup\Lambda_-.
\]
The preceding comparison gives $\Lambda_-\ne\varnothing$, and the definition shows that $\Lambda_-$ is downward closed in $(-\infty,0)$.

\begin{lemma}\label{Lem4.4}
At the limiting position,
\[
u\le u_{\lambda_0} \quad\text{in }\Sigma_{\lambda_0}\setminus\Gamma_{\lambda_0}.
\]
\end{lemma}

\begin{proof}
Choose $\lambda_j\in\Lambda_-$ with $\lambda_j\uparrow\lambda_0$. If $x\in\Sigma_{\lambda_0}\setminus\Gamma_{\lambda_0}$, then $x\in\Sigma_{\lambda_j}\setminus\Gamma_{\lambda_j}$ for all large $j$, because $\Gamma_{\lambda_j}$ converges to $\Gamma_{\lambda_0}$ in the Hausdorff metric. Hence
\[
u(x)\le u(x^{\lambda_j}).
\]
The assertion follows by continuity as the reflected points converge to $x^{\lambda_0}\notin\Gamma$.
\end{proof}

\begin{lemma}\label{Lem4.5}
Assume $\lambda_0<0$. There are $R>1$, $\delta\in(0,-\lambda_0/2)$, $\rho>0$, and $\tau\in(0,\min\{\delta/2,\rho/4\})$ such that, for every $\lambda\in(\lambda_0,\lambda_0+\tau)$,
\begin{equation}\label{eq4.1}
E_\lambda\subset \left(B_R\cap\{|x_1-\lambda_0|<\delta\}\right) \cup \left(B_R\cap O_\rho(\Gamma_{\lambda_0})\right) \cup (\R^N\setminus B_R),
\end{equation}
and
\[
\mathfrak M_u(E_\lambda)<\kappa_*,
\]
where $\kappa_*$ is the constant in Proposition~\ref{Prop3.5}.
\end{lemma}

\begin{proof}
By Lemma~\ref{Lem4.4}, Proposition~\ref{Prop4.2}, and Lemma~\ref{Lem4.3},
\begin{equation}\label{eq4.2}
u<u_{\lambda_0} \quad\text{in }\Sigma_{\lambda_0}\setminus\Gamma_{\lambda_0}.
\end{equation}
Choose $\varepsilon>0$ so small that
\[
(3\varepsilon)^{2\hstar-2}+3\varepsilon<\kappa_*.
\]
By \eqref{eq2.6}, choose $R$ so large that
\[
\|u\|_{L^{\sstar}(\R^N\setminus B_R)}<\varepsilon, \qquad \|\cA_u\|_{L^{N/(2s)}(\R^N\setminus B_R)}<\varepsilon.
\]
The sets $\Gamma$ and $\Gamma_{\lambda_0}$ have distance $2|\lambda_0|$. Hence, for sufficiently small $\delta$ and $\rho$, both
\[
B_R\cap\{|x_1-\lambda_0|<\delta\} \quad\text{and}\quad B_R\cap O_\rho(\Gamma_{\lambda_0})
\]
lie in a fixed compact subset of $\R^N\setminus\Gamma$. Their measures tend to zero as $\delta,\rho\downarrow0$, since $|\Gamma_{\lambda_0}|=0$. Absolute continuity of the two critical integrals allows us to choose $\delta$ and $\rho$ such that the $L^{\sstar}$-norm of $u$ and the $L^{N/(2s)}$-norm of $\cA_u$ on each set are less than $\varepsilon$.

Set
\[
K=\overline{B_R}\cap\{x_1\le\lambda_0-\delta\} \setminus O_\rho(\Gamma_{\lambda_0}).
\]
If $K$ is nonempty, then $K\Subset\Sigma_{\lambda_0}\setminus\Gamma_{\lambda_0}$. By \eqref{eq4.2} and continuity,
\[
m=\min_K(u_{\lambda_0}-u)>0.
\]
For $x\in K$ and $0<\lambda-\lambda_0<\rho/4$,
\[
\dist(x^\lambda,\Gamma)\ge \dist(x^{\lambda_0},\Gamma)-|x^\lambda-x^{\lambda_0}|\ge\rho-2(\lambda-\lambda_0)>\rho/2.
\]
For $0<\tau<\rho/4$, the set $\bigcup_{\lambda\in[\lambda_0,\lambda_0+\tau]}K^\lambda$ is compact and disjoint from $\Gamma$. Uniform continuity of $u$ on this set gives
\[
u<u_\lambda\quad\text{on }K
\]
for $\lambda\in(\lambda_0,\lambda_0+\tau)$, after choosing $0<\tau<\min\{\delta/2,\rho/4\}$ sufficiently small. If $K$ is empty, take any $\tau$ in this interval.

Let $x\in E_\lambda\cap B_R$ and suppose that $x\notin O_\rho(\Gamma_{\lambda_0})$. Then $x\notin K$, so $x_1>\lambda_0-\delta$. Since $x_1<\lambda<\lambda_0+\tau<\lambda_0+\delta$, one has $|x_1-\lambda_0|<\delta$. This proves \eqref{eq4.1}. The choices above give
\[
\|u\|_{L^{\sstar}(E_\lambda)}<3\varepsilon, \qquad \|\cA_u\|_{L^{N/(2s)}(E_\lambda)}<3\varepsilon,
\]
and therefore $\mathfrak M_u(E_\lambda)<\kappa_*$.
\end{proof}

\begin{proposition}\label{Prop4.6}
One has $\lambda_0=0$. Consequently,
\begin{equation}\label{eq4.3}
u(x)\le u_\lambda(x) \quad\text{for every }\lambda<0,\quad x\in\Sigma_\lambda\setminus\Gamma_\lambda.
\end{equation}
\end{proposition}

\begin{proof}
Assume $\lambda_0<0$. Lemma~\ref{Lem4.5} and Proposition~\ref{Prop3.5} give
\[
u\le u_t \quad\text{in }\Sigma_t\setminus\Gamma_t \quad\text{for every }t\in(\lambda_0,\lambda_0+\tau).
\]
For every $t<\lambda_0$, choose $\bar\lambda\in\Lambda_-$ with $t<\bar\lambda<\lambda_0$. The comparison at $t$ is part of the defining condition for $\bar\lambda\in\Lambda_-$. At $t=\lambda_0$, the comparison follows from Lemma~\ref{Lem4.4}. Therefore, for each $\lambda\in(\lambda_0,\lambda_0+\tau)$, the comparison holds for every $t\le\lambda$, so $\lambda\in\Lambda_-$. This contradicts the definition of $\lambda_0$. Thus $\lambda_0=0$. Finally, for any $\lambda<0$, the equality $\sup\Lambda_-=0$ provides $\bar\lambda\in\Lambda_-$ with $\lambda<\bar\lambda<0$. Thus $\lambda\in\Lambda_-$ by definition, proving \eqref{eq4.3}.
\end{proof}

\begin{proposition}\label{Prop4.7}
For every $\lambda>0$,
\begin{equation}\label{eq4.4}
u(x)\le u_\lambda(x) \quad\text{for }x_1>\lambda, \quad x\notin\Gamma_\lambda.
\end{equation}
\end{proposition}

\begin{proof}
Apply Propositions~\ref{Prop3.5} and~\ref{Prop4.6} to the function
\[
\widetilde u(x_1,x')=u(-x_1,x').
\]
The equation, the two tail conditions, the capacity condition, and the existence of a nonremovable singular point are invariant under this reflection. Rewriting the comparison for $\widetilde u$ in the left half-space gives \eqref{eq4.4}.
\end{proof}

\begin{proposition}\label{Prop4.8}
Let $u$ be a positive weak singular solution of \eqref{eq1.1}. Assume $\Gamma\subset\{x_1=0\}$ is compact with zero $(s,2)$-capacity, $u$ has a nonremovable singular point on $\Gamma$, and the two tail conditions \eqref{eq2.5} hold. Then
\[
u(x_1,x')=u(-x_1,x') \quad\text{in }\R^N\setminus\Gamma,
\]
and $u$ is strictly increasing in $x_1$ on $\{x_1<0\}$ and strictly decreasing in $x_1$ on $\{x_1>0\}$.
\end{proposition}

\begin{proof}
Let $x=(x_1,x')$ with $x_1<0$, and set $\bar x=(-x_1,x')$. Choose $\lambda_j\uparrow0$ with $x_1<\lambda_j<0$ and $x\notin\Gamma_{\lambda_j}$. Proposition~\ref{Prop4.6} gives
\[
u(x)\le u(x^{\lambda_j}).
\]
Passing to the limit yields $u(x)\le u(\bar x)$. Applying Proposition~\ref{Prop4.7} to $\bar x$ and taking planes decreasing to zero gives the reverse inequality. Hence $u$ is symmetric with respect to $\{x_1=0\}$.

Let $x_1<y_1<0$ and put $\lambda=(x_1+y_1)/2$. Then $(x_1,x')\in\Sigma_\lambda\setminus\Gamma_\lambda$ and its reflection is $(y_1,x')$. The weak inequality follows from Proposition~\ref{Prop4.6}; it is strict by Proposition~\ref{Prop4.2} and Lemma~\ref{Lem4.3}. Strict decrease on the right follows from symmetry.
\end{proof}

\section{Proof of the Main Result}\label{Sec5}

We now deduce symmetry and monotonicity of the original solution from the corresponding properties of its Kelvin transforms. Let
\[
\mathscr R(x_1,x')=(-x_1,x')
\]
be reflection across $\{x_1=0\}$. If $a\in\{x_1=0\}$, then
\begin{equation}\label{eq5.1}
\mathcal I_a\circ\mathscr R=\mathscr R\circ \mathcal I_a, \qquad |\mathscr R x-a|=|x-a|.
\end{equation}

\begin{lemma}\label{Lem5.1}
Let $a\in\{x_1=0\}\setminus\Gamma$. If the Kelvin transform $u_a$ is symmetric with respect to $\{x_1=0\}$, then $u$ has the same symmetry.
\end{lemma}

\begin{proof}
If $y=a$, the assertion is immediate because $a$ lies on the reflecting hyperplane. Let $y\in\R^N\setminus(\Gamma\cup\{a\})$ and put $x=\mathcal I_a(y)$. Since $\mathcal I_a$ is an involution,
\[
u_a(x)=|x-a|^{-(N-2s)}u(y) =|y-a|^{N-2s}u(y).
\]
Using \eqref{eq5.1},
\[
u_a(\mathscr R x)=|y-a|^{N-2s}u(\mathscr R y).
\]
By symmetry of $u_a$, the left-hand sides are equal. Cancelling the common factor on the right proves $u(y)=u(\mathscr R y)$.
\end{proof}

\begin{lemma}\label{Lem5.2}
Assume that, for every $a\in\{x_1=0\}\setminus\Gamma$, the Kelvin transform $u_a$ is nondecreasing in the first coordinate on $\{x_1<0\}$. Then
\begin{equation}\label{eq5.2}
u(r,x')\le u(t,x') \quad\text{whenever }r<t<0.
\end{equation}
\end{lemma}

\begin{proof}
Fix $r<t<0$ and $x'\in\R^{N-1}$. Choose a unit vector $e\in\R^{N-1}$ and, for sufficiently large integers $k$, set
\[
a_k=(0,x'+ke),\qquad D_k=k^2+r^2,\qquad F_k=k^2+t^2.
\]
Since $\Gamma$ is compact, $a_k\notin\Gamma$ for large $k$. Define
\[
X_k=\left(\frac r{D_k},\,x'+ke-\frac{ke}{D_k}\right), \qquad Y_k=\left(\frac t{D_k},\,x'+ke-\frac{ke}{D_k}\right).
\]
Then $X_k=\mathcal I_{a_k}(r,x')$. The two points have the same tangential coordinates and satisfy $(X_k)_1<(Y_k)_1<0$, so
\[
u_{a_k}(X_k)\le u_{a_k}(Y_k).
\]
A direct calculation gives
\[
Q_k:=\mathcal I_{a_k}(Y_k) =\left(\frac{tD_k}{F_k},\,x'+\frac{k(t^2-r^2)}{F_k}e\right) \longrightarrow(t,x').
\]
Moreover,
\[
|X_k-a_k|^2=D_k^{-1},\qquad |Y_k-a_k|^2=F_kD_k^{-2}.
\]
The definition of the Kelvin transform therefore yields
\[
u(r,x')\le\left(\frac{D_k}{F_k}\right)^{\frac{N-2s}{2}}u(Q_k).
\]
Since $D_k/F_k\to1$ and $(t,x')\notin\Gamma$, continuity of $u$ proves \eqref{eq5.2}. The argument uses monotonicity for each fixed Kelvin center and requires no uniform estimate as the centers tend to infinity.
\end{proof}

\begin{proof}[Proof of Theorem~\ref{Thm1.2}]
Choose any $a\in\{x_1=0\}\setminus\Gamma$. Proposition~\ref{Prop2.3} shows that the Kelvin transform $u_a$ is a positive weak singular solution of the same equation, with a compact singular set of zero capacity
\[
\Gamma_a=\mathcal I_a(\Gamma)\cup\{a\}\subset\{x_1=0\}.
\]
It has a nonremovable singular point. Proposition~\ref{Prop2.6}(ii) supplies both tail conditions in \eqref{eq2.5}, without requiring either condition for $u$. Proposition~\ref{Prop4.8} therefore applies to $u_a$ for every $0<\mu<N$. Lemma~\ref{Lem5.1} then shows that $u$ has the same hyperplane symmetry as $u_a$.

Since $a$ was arbitrary, Propositions~\ref{Prop2.3}, \ref{Prop2.6}(ii), and~\ref{Prop4.8} show that the Kelvin transform $u_a$ about every regular center is strictly increasing in the first coordinate on $\{x_1<0\}$. Lemma~\ref{Lem5.2} therefore gives weak transverse monotonicity on the negative half-space:
\[
u(r,x')\le u(t,x') \quad\text{for }r<t<0.
\]
To prove the reflected comparison for each plane $T_\lambda$, $\lambda<0$, we combine this monotonicity with symmetry. Let $x_1<\lambda$ and put $y_1=2\lambda-x_1$. If $y_1<0$, we use the preceding monotonicity directly:
\[
u(x_1,x')\le u(y_1,x').
\]
If $y_1=0$ and $(0,x')\notin\Gamma$, the same conclusion follows by taking a limit from the negative side. If $y_1=0$ and $(0,x')\in\Gamma$, then $(x_1,x')\in\Gamma_\lambda$, so this point is excluded from the comparison domain. If $y_1>0$, symmetry gives
\[
\begin{gathered}
u(y_1,x')=u(-y_1,x'),\\
x_1<-y_1=x_1-2\lambda<0.
\end{gathered}
\]
Thus weak monotonicity again yields
\[
u(x_1,x')\le u(y_1,x').
\]
Therefore,
\begin{equation}\label{eq5.3}
u\le u_\lambda \quad\text{in }\Sigma_\lambda\setminus\Gamma_\lambda \quad\text{for every }\lambda<0.
\end{equation}

Proposition~\ref{Prop4.2} does not require \eqref{eq2.5} and therefore applies to the original solution. Lemma~\ref{Lem4.3} excludes the case $u\equiv u_\lambda$ because $u$ has a nonremovable singular point. Hence \eqref{eq5.3} is strict in $\Sigma_\lambda\setminus\Gamma_\lambda$ for every $\lambda<0$. Taking
\[
\lambda=\frac{x_1+y_1}{2}
\]
gives
\[
u(x_1,x')<u(y_1,x') \quad\text{whenever }x_1<y_1<0.
\]
Strict decrease in the positive half-space follows from the already proved symmetry.
\end{proof}

\begin{proof}[Proof of Corollary~\ref{Cor1.3}]
Let $e\in\R^N$ be a unit vector. Since $\Gamma=\{0\}\subset e^\perp$ and the equation and all the assumptions of Theorem~\ref{Thm1.2} are invariant under rotations, that theorem applies in the direction $e$. Hence $u$ is symmetric with respect to every hyperplane through the origin, and therefore $u(x)=U(|x|)$ for some function $U:(0,\infty)\to(0,\infty)$. For $0<r<t$, we have $U(t)=u(t,0)<u(r,0)=U(r)$ by the strict monotonicity in Theorem~\ref{Thm1.2}, with $x'=0$.
\end{proof}

\section*{Acknowledgments}
This work is supported by the National Natural Science Foundation of China (12301145, 12561020, 12261107) and Yunnan Fundamental Research Projects (202401AU070123, 202601AT070048).
The authors acknowledge the use of ChatGPT (OpenAI) for exploratory discussions and language refinement during the preparation of this work. All mathematical statements, proofs, and references were independently checked by the authors, who take full responsibility for the content of the paper.

\medskip
{\bf Data availability:} Data sharing is not applicable to this article as no new data were created or analyzed in this study.

\medskip
{\bf Conflict of interest:} The authors declare that they have no conflict of interest.

\end{document}